\documentclass[12pt, a4paper]{myarticle}

\usepackage{my_proof}

\usepackage{ifpdf}
\ifpdf
 \usepackage[pdftex]{graphicx}
\else
 \usepackage[dvips]{graphicx}
 \DeclareGraphicsExtensions{.png, .gif , .jpg, .eps}
\fi

\usepackage{epic}

\usepackage{amsmath}                                                            
\usepackage{amssymb}                                                            
\usepackage{amsfonts}                                                           

\usepackage{epsfig}                                                             
                                                                                
\usepackage{ifthen}                                                             

\newcounter{bsplemsatz}
\renewcommand{\thebsplemsatz}{\arabic{bsplemsatz}}

\newenvironment{env}[2][]{%
  \refstepcounter{bsplemsatz}                                                   
  \trivlist                                                                     
  \item[\hskip\labelsep{\bf #2\relax\ \thebsplemsatz.}%
  ]\ignorespaces\itshape%
}{%
  \endtrivlist                                                                  
}                                                                               

\newcommand{\RR}{\mathbb{R}}

\DeclareMathAlphabet{\masc}{U}{eus}{m}{n}                                       
\DeclareMathAlphabet{\mafr}{U}{euf}{m}{n}                                       
\DeclareMathAlphabet{\malf}{OT1}{cmtt}{m}{it}               

   \def\soft#1{\leavevmode\setbox0=\hbox{h}\dimen7=\ht0\advance                                                         
    \dimen7 by-1ex\relax\if t#1\relax\rlap{\raise.6\dimen7                                                              
    \hbox{\kern.3ex\char'47}}#1\relax\else\if T#1\relax                                                                 
    \rlap{\raise.5\dimen7\hbox{\kern1.3ex\char'47}}#1\relax                                                             
    \else\if d#1\relax\rlap{\raise.5\dimen7\hbox{\kern.9ex                                                              
    \char'47}}#1\relax\else\if D#1\relax\rlap{\raise.5\dimen7                                                           
    \hbox{\kern1.4ex\char'47}}#1\relax\else\if l#1\relax                                                                
    \rlap{\raise.5\dimen7\hbox{\kern.4ex\char'47}}#1\relax                                                              
    \else\if L#1\relax\rlap{\raise.5\dimen7\hbox{\kern.7ex                                                              
    \char'47}}#1\relax\else\message{accent \string\soft                                                                 
    \space #1 not defined!}#1\relax\fi\fi\fi\fi\fi\fi}                                        
    
\newcommand{\obda}{w.l.o.g.~}

\newcommand{\cclass}{{\malf C}}

\newcommand{\patterngr}[1]{\Gamma(#1)}
\newcommand{\support}[2]{\operatorname*{supp}(#1;#2)}

\begin{document}

\title{Tree Decomposition By Eigenvectors}
\author{T. Sander\and J. W. Sander}

\date{\today\\[1em]
\emph{Institut f\"ur Mathematik\\
Technische Universit\"at Clausthal\\
D-38678 Clausthal-Zellerfeld, Germany}\\[1em]
{e-mail:} \{torsten$\vert$juergen\}.sander@math.tu-clausthal.de}

\maketitle



\thispagestyle{empty}

\centerline{\large \bf Abstract}
In this work a composition-decomposition technique is presented that correlates tree eigenvectors with
certain eigenvectors of an associated so-called skeleton forest. In particular, the matching
properties of a skeleton determine the multiplicity of the corresponding tree eigenvalue. 
As an application a characterization of trees that admit eigenspace bases with entries only 
from the set $\{0,1,-1\}$ is presented. Moreover, a result due to Nylen concerned with partitioning
eigenvectors of tree pattern matrices is generalized.






\section{Introduction}

Eigenspaces of graphs have been researched to some degree since many years. This is especially the case for
the null space, which has been studied for a number of graph classes. 
Compared to the amount of research spent on the spectrum of graphs, only little 
attention has been given to what the eigenvectors of graphs really look like. 
There exist explicit results on paths, cycles, circulant graphs, hypercubes and some other graph classes.

This work investigates eigenvectors of trees. The eigenvectors of a graph, i.e.~eigenvectors
of the adjacency matrix, can be considered as real valued functions on its vertex set. 
One may partition the vertices of a tree by grouping those vertices on which there exist non-zero values for
some vector from a fixed eigenspace and those on which every vector from that space vanishes.
Then the components arising from vertices of the first kind can be contracted into single vertices
(these components are loosely related to the so-called nut graphs studied in \cite{sciriha98} and \cite{scigut98}). Together
with any adjacent vertices they induce a so-called skeleton forest, a concept initially hinted at in \cite{nylen98}. 
We show that the null space of this skeleton provides a blue print for the vectors from the considered eigenspace of the original tree. 
Moreover, its matching properties can be utilized to determine the eigenspace dimension. 
A converse technique allows to create trees with certain eigenspace properties from given skeletons.

The partitioning of graphs according to eigenvector structure has been studied before, but with different aims.
For example in \cite{powers88} the components induced by the signs of the entries of eigenvectors 
have been investigated and an upper bound for the number of zero entries in eigenvectors of graphs 
has been derived.

Let us point out some more related research. 
An overview of theorems on the null space dimension of trees
and bipartite graphs in general can be found in \cite{cvetgut72}. A construction
of trees with a given number of vertices and upper vertex degree bound that
achieve the greatest possible null space dimension is given in \cite{fiogutsci05}.
For line graphs of trees it is known \cite{sciriha99} that 
the order of eigenvalue zero can only be either 0 or 1. General research on eigenvectors of trees
has been conducted, for example, in \cite{bevisdom95}, \cite{fricke96} and \cite{johnschild96}.
In particular, it is shown that the number of vertices covered by a maximum matching of a tree
equals the rank of its adjacency matrix. The authors also present and analyze efficient algorithms for determining
rank, characteristic polynomial and eigenvectors of trees. 

The case of finding particularly simple eigenspace bases is a rather new research topic. Owing to
sparsity, bases that only contain entries from the set $\{0, 1, -1\}$ are of particular interest.
Such bases we call simply structured. 
For the null space of the incidence matrix of a tree the situation is comparably simple.
It has been shown in \cite{akbarigha06} that every graph without cut-edge has a simply structured incidence matrix null space basis.
Further, there exists a characterization of all graphs with a simply structured incidence matrix null space basis. 
For the adjacency matrix of a tree so far it has been independently shown in \cite{akbariali06}, \cite{fox1} that every tree has such a simply structured null space basis. 
The only other two feasible eigenvalues that allow the construction of such bases are $1$ and $-1$.
As an application of our decomposition technique, we settle the question of characterizing all trees whose eigenspaces admit simply 
structured bases. 

As a second application we rediscover and generalize a result published in \cite{nylen98}.
Associate with a given real $n\times n$ symmetric matrix an undirected graph $\Gamma(A)$ on $n$ vertices $1,2,\ldots,n$ by including the edge
connecting vertex $i$ to vertex $j$ in the edge set if and only if the entry at position $(i,j)$ in $A$ is non-zero.
If $\Gamma(A)$ is a tree, then the dimension of the null space of $A$ is the number of connected components of the subgraph of $\Gamma(A)$ induced by 
the set of indices $i$ such that $x_i$ is non-zero for some vectors $x$ from the null space $A$,
minus the number of indices adjacent to (with respect to $\Gamma(A)$) but not in that index set.
Here the topic of acyclic matrices is touched. A number of results exist on the spectra of such matrices in general. 
One of the first papers in that area was \cite{fiedler75}, where it is shown that the signs of eigenvector entries are related to the position of 
the corresponding eigenvalue in the ordered spectrum of a tree pattern matrix.
In \cite{johnsaia02} it is shown that it is in general not possible to express
maximum eigenvalue multiplicities of a class of acyclic matrices given by some tree
in terms of the degrees of the vertices. It is moreover well known that
any tree $T$ has at least diam(T)+1 distinct eigenvalues. This fact can be generalized
to tree pattern matrices \cite{lealjohn02}. For even more general results that
cover conditions on the solution of linear equations or
existence of certain eigenvalues and eigenvectors for a class of
matrices whose zero-nonzero pattern matches a given digraph see e.g.~\cite{macdonald93}.

\section{Basics and Notation}

In this paper we will only consider finite, loopless, simple graphs. The basics of graph theory
are treated e.g.~in \cite{diestel}, \cite{harary}. We will only introduce some initial notation
and will state further definitions along the way.

Let $G$ be a graph. The vertex set of $G$ will be written as $V(G)$. 
Given a set $M\subseteq V(G)$ we denote by $G\setminus M$ the graph formed by
removing the vertices of $M$ and all their adjacent edges from $G$.

When we talk of a {\em vertex bipartition} of a bipartite graph (e.g.~a forest) we mean a disjoint partition
of its vertex set into two sets such that every edge of the graph runs from a vertex of the first set
to a vertex of the second set.

The foundations of algebraic graph theory can be found in \cite{biggs93}, \cite{eigenspacesofgraphs}, \cite{spectraofgraphs}, \cite{godsil01}. 
The eigenvalues of a graph $G$ with vertex set $V=\{v_1,\ldots,v_n\}$ are the eigenvalues of its adjacency matrix $A=(a_{ij})$
which is defined by $a_{ij}=1$ if $v_i$ is adjacent to $v_j$ and $a_{ij}=0$ otherwise. Note that this eigenvalue definition
is independent of the chosen vertex order. Eigenvalues of graphs are real. For a bipartite graph $\lambda$ is
an eigenvalue if and only if $-\lambda$ is an eigenvalue of the graph as well.

Suppose that $Ax=\lambda x$, where $x=(x_1,\ldots,x_n)^T$. 
If we assign value $x_i$ to vertex $v_i$, it is easily seen that  for every vertex  the sum over the
values of its neighbors equals $\lambda$ times its own value. We will hereafter refer to this as the {\em summation rule}.

Note that we do not consider the null vector an eigenvector although it formally belongs to an eigenspace.

We conclude this section by quoting a basic result that we will frequently refer to:

\begin{env}{Lemma}\cite{fiedler75}\label{lemfiedler}
Let $T$ be a tree. Let $v$ be an eigenvector of $T$ for eigenvalue $\lambda$. If $v$ does not have
any zero entries, then $\lambda$ necessarily has multiplicity one.
\end{env}

\section{Main Results}\label{chmain}

\subsection{Tree eigenvector decomposition}

\newcommand{\comp}[1]{{\mafr C}(#1)}
\newcommand{\nzero}[2]{{N_{#2}(#1)}}
\newcommand{\nzerox}[2]{{N(#1,#2)}}
\newcommand{\nzeroc}[2]{{N^C_{#2}(#1)}} 
\newcommand{\nzerocx}[2]{{N^C(#1,#2)}}
\newcommand{\skelx}[2]{S(#1,#2)}
\newcommand{\skel}[2]{S_{#2}(#1)}

Let $G$ be a graph and $M=\{X_1,\ldots,X_r\}$ a set of mutually vertex disjoint subgraphs of $G$. Then by $G/\{X_i\}_{i=1}^r$ or $G/M$ we denote 
the graph that results from the contraction of each subgraph $X_i$ in $G$ to a single vertex $x_i$. Further, let 
$\comp{G}$ denote the set of components of $G$.

Let $x$ be an eigenvector for eigenvalue $\lambda$ of graph $G$. Let $\nzerox{G}{x}$ be the set of those vertices 
of $G$ on which $x$ vanishes. 
Moreover, let $\nzero{G}{\lambda}$ mean the set of vertices on which every non-zero eigenvector 
for eigenvalue $\lambda$ of $G$ vanishes.

\begin{env}{Lemma}\label{restrictev}
Let $G$ be a graph and $x$ an eigenvector for its eigenvalue $\lambda$. Then:
\begin{enumerate}
\item For any $C\in\comp{G\setminus \nzerox{G}{x}}$ the restriction $x\vert_{C}$ is an eigenvector of the graph $C$ 
for eigenvalue $\lambda$. If $G$ is a tree, then 
$x\vert_{C}$ constitutes an eigenspace basis of the subtree $C$ for eigenvalue $\lambda$.
\item For any $C\in\comp{G\setminus \nzero{G}{\lambda}}$ the restriction $x\vert_{C}$ is either the null vector or an eigenvector of the graph $C$ 
for eigenvalue $\lambda$. If $G$ is a tree and $x\vert_{C}\not=0$, then 
$x\vert_{C}$ does not contain any zero entries and, moreover, constitutes an eigenspace basis of the subtree $C$ for eigenvalue $\lambda$.
\end{enumerate}
\end{env}

\begin{proof}
In the first case the claim follows directly from the definition of $\nzerox{G}{x}$, the summation rule
and Lemma \ref{lemfiedler}.
The case $C\in\comp{G\setminus \nzero{G}{\lambda}}$ is similar, with only one additional argument. 
Let $v_1,\ldots,v_k$ be the vertices of $C$. For every vertex $v_i$ there exists an
eigenvector $x_i$ of $T$ for eigenvalue $\lambda$ whose restriction $x_i\vert_{C}$ does not vanish on $v_i$.
It is straightforward to show (see e.g.~Lemma~7 in \cite{nylen98}) 
that there exists a linear combination of these vectors $x_i$
that has no zero entries so that by Lemma \ref{lemfiedler} the associated 
eigenvalue $\lambda$ of $C$ has multiplicity one.
\end{proof}

It is easy to see that the eigenspace basis related claim of the lemma does not extend to general graphs. For example, look at the complete graph $K_3$ and its
eigenvectors $(2,-1,-1)$, $(-1,2,-1)$ for eigenvalue $-1$. 

Let $x$ be an eigenvector for eigenvalue $\lambda$ of a given tree $T$. 
Further let $C_i$, $i=1,\ldots,r$, be the elements of $\comp{T\setminus \nzerox{T}{x}}$.
We will now concentrate on a particularly interesting subset of $\nzerox{T}{x}$. Namely, 
let $\nzerocx{T}{x}$ consist of all those vertices of $\nzerox{T}{x}$ that are adjacent to at least one of the subgraphs $C_i$ in $T$.

\begin{env}{Lemma}\label{skelxtree}
Let $T$ be a tree and $x$ an eigenvector for its eigenvalue $\lambda$.
Let $C_i$, $i=1,\ldots,r$, be the elements of $\comp{T\setminus \nzerox{T}{x}}$.
Further let $c_i$ denote the associated contracted vertices of $T/\{C_i\}_{i=1}^r$.

Then the vertex set $\nzerocx{T}{x} \cup \{c_1,\ldots,c_r\}$ induces a forest $F$ in $T/\{C_i\}_{i=1}^r$ such that
the leaves of $F$ form a subset of $\{c_1,\ldots,c_r\}$ and are also leaves of $T/\{C_i\}_{i=1}^r$.
\end{env}

\begin{proof}
Clearly, the contraction $T/\{C_i\}_{i=1}^r$ of the tree $T$ by the sub-forest $\bigcup\{C_i\}$ is a forest. 
So the induced subgraph $F$ of $T/\{C_i\}_{i=1}^r$ must be a forest as well.

Now consider an element $v$ of $\nzerocx{T}{x}$. By construction and since $T$ is a tree there exists a one-to-one mapping of 
the non-zero weight neighbors of $v$ to a subset of $\comp{T\setminus \nzerox{T}{x}}$.
By definition, $v$ is adjacent to at least one component $C_i$, but since the sum over the neighbors of $v$ must 
vanish we see that it must be adjacent to at least two such components. Consequently, $v$ is adjacent to
at least two of the vertices $c_i$ in both $T/\{C_i\}_{i=1}^r$ and $F$. So the leaves in $F$
are a subset of $\{c_1,\ldots,c_r\}$.

Assume that $c_k$ is a leaf of $F$ that is not a leaf of $T/\{C_i\}_{i=1}^r$. Then in $T/\{C_i\}_{i=1}^r$ there would exist
a neighbor $w$ of $c_k$ such that $w\in\nzerox{T}{x}\setminus\nzerocx{T}{x}$. But then $w$ could only be
adjacent to zero-weight vertices, a contradiction.
\end{proof}

In \cite{sciriha98}, \cite{scigut98} graphs with nullity one and corresponding eigenvector without zero entries,
so-called nut graphs, are studied. Nut graphs have a number of interesting properties. 
Clearly, the components of $T\setminus \nzerox{T}{x}$ are nut graphs if $x$ is an eigenvector for
eigenvalue $0$. However, the theory on nut graphs does not yield any insight in the case of trees because
it is easy to see that for a tree $K_1$ is the only possible nut graph.

In the following, let $\skelx{T}{x}$ denote the forest $F$ by Lemma \ref{skelxtree} associated with a given
tree $T$ and eigenvector $x$. We call $\skelx{T}{x}$ the {\em $x$-skeleton} of $T$.

\begin{figure}[ht]
\begin{center}
\epsfig{file=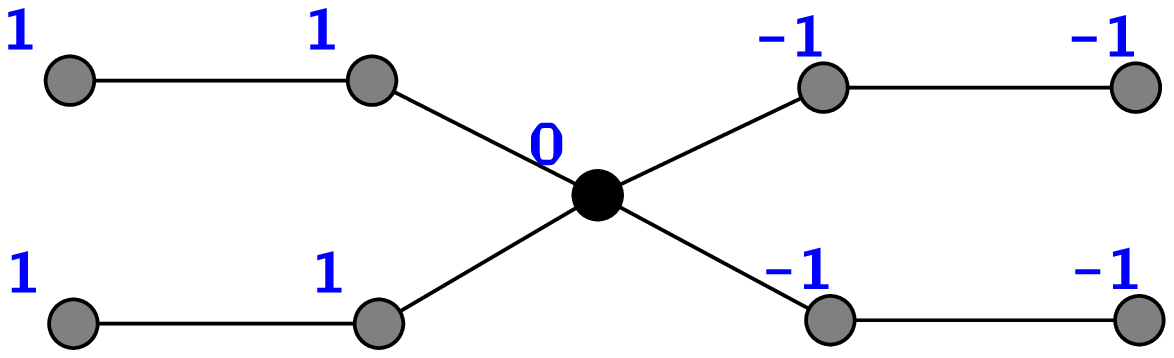,scale=0.35}\hspace*{20mm}
\epsfig{file=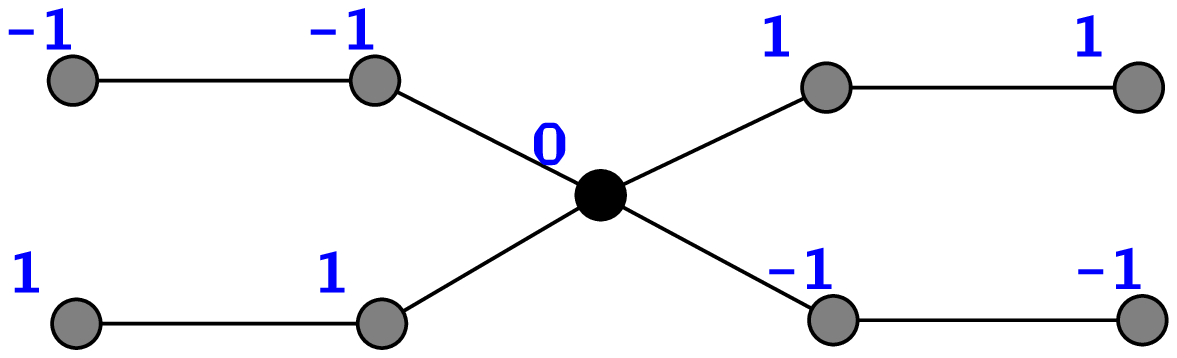,scale=0.35}
\end{center}
\caption{Eigenvectors with the same $x$-skeleton} \label{figsameskel}
\end{figure}

Note that $x$-skeletons do not characterize an eigenspace basis, i.e.~there may exist
linearly independent eigenvectors $x,x'$ for eigenvalue $\lambda$ of a tree $T$ such that
$\comp{T\setminus \nzerox{T}{x}} = \comp{T\setminus \nzerox{T}{x'}}$. An example is shown in Figure \ref{figsameskel}.

\begin{env}{Theorem}\label{nzerocompcontain}
Let $T$ be a tree and $x$ an eigenvector for eigenvalue $\lambda$ of $T$. Then,
\[
  \comp{T\setminus \nzerox{T}{x}} \subseteq \comp{T\setminus \nzero{T}{\lambda}}.
\]
\end{env}

\begin{proof}
Let $C\in\comp{T\setminus \nzerox{T}{x}}$.
Clearly, none of the vertices of $C$ belong to the set $\nzero{T}{\lambda}$. Therefore
$C$ is a subgraph of some component $C'\in\comp{T\setminus \nzero{T}{\lambda}}$.

According to Lemma \ref{restrictev} the vectors $x\vert_{C}$ and $x\vert_{C'}$ are eigenvectors for eigenvalue $\lambda$
on $C$ and $C'$, respectively, and in particular do not have any zero entries on $C'$.
Hence $C=C'$.
\end{proof}

\begin{env}{Corollary}\label{nonoverlapci}
Let $x,x'$ be linearly independent eigenvectors for eigenvalue $\lambda$ of a given tree and let
$C\in\comp{T\setminus \nzerox{T}{x}}$, $C'\in\comp{T\setminus \nzerox{T}{x'}}$. Then either $C$ and $C'$ are identical or they 
are disjoint subgraphs of $T$.
\end{env}

\begin{env}{Corollary}\label{nzerocxinnzerox}
Let $T$ be a tree with eigenvector $x$ for eigenvalue $\lambda$. Then,
\[ 
  \nzerocx{T}{x}\subseteq\nzero{T}{\lambda}.
\]
\end{env}

\begin{env}{Corollary}\label{compunion}
Let $T$ be a tree with eigenvalue $\lambda$. Then,
\[
   \comp{T\setminus \nzero{T}{\lambda}} = \bigcup\limits_x \comp{T\setminus \nzerox{T}{x}},
\]
where the union is taken over all eigenvectors $x$ for eigenvalue $\lambda$ of $T$.
\end{env}

As a consequence of Corollary \ref{compunion}
we can safely merge the $x$-skeleton forests of an entire eigenspace.
Let $T$ be a tree and let $C_1,\ldots,C_r$ be the elements of $\comp{T\setminus \nzero{T}{\lambda}}$.
Let the associated contracted vertices
in $T/\{C_i\}_{i=1}^r$ be $c_1,\ldots,c_r$. Denote the union of the sets $\nzerocx{T}{x}$ by $\nzeroc{T}{\lambda}$.
Now we define the {\em skeleton} $\skel{T}{\lambda}$ as the sub-forest of $T/\{C_i\}_{i=1}^r$ induced by the vertices of
$\nzeroc{T}{\lambda} \cup \{c_1,\ldots,c_r\}$.

\begin{figure}[ht]
\begin{center}
\epsfig{file=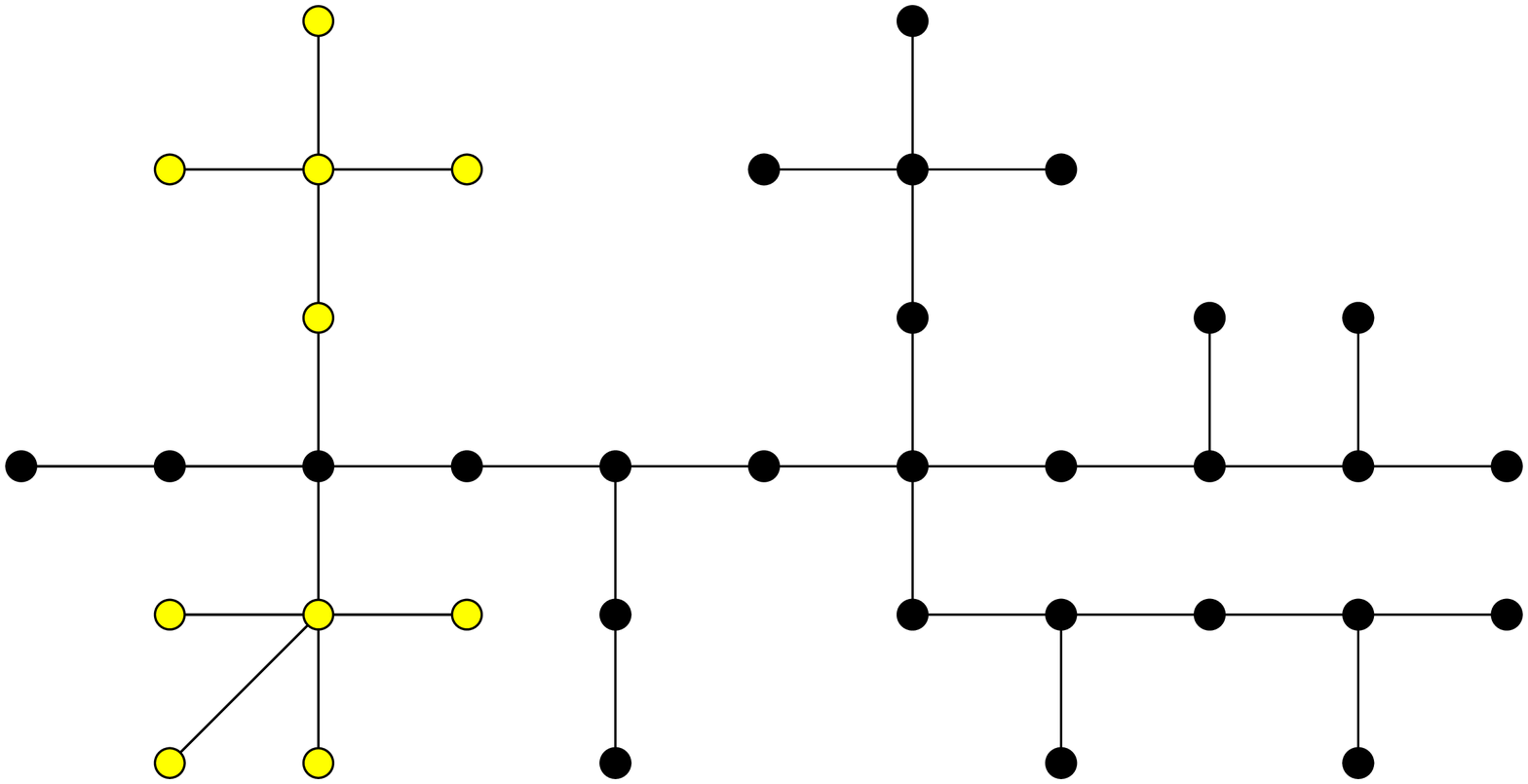,scale=0.25}\hspace*{3mm}
\epsfig{file=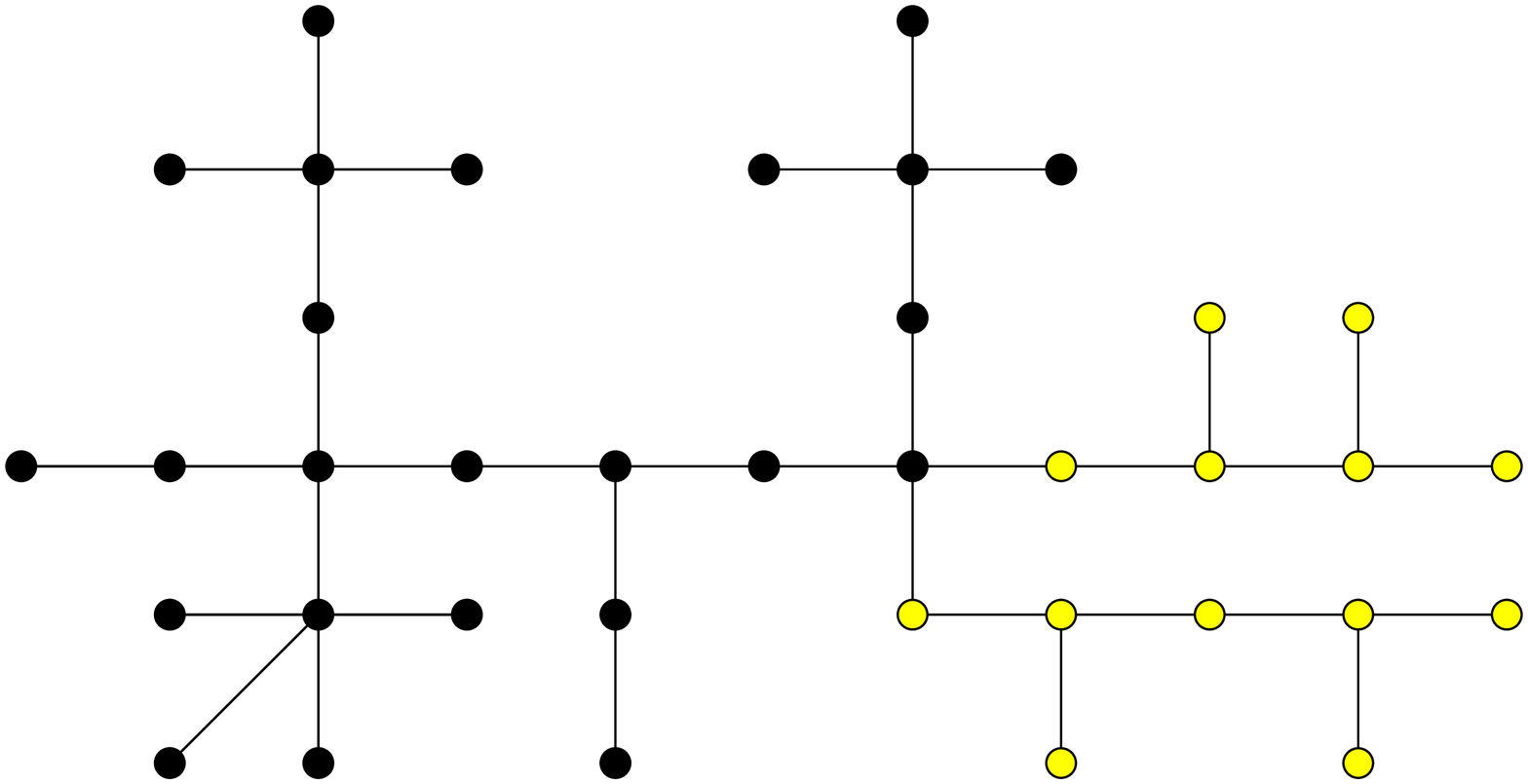,scale=0.25}\\[5mm]
\epsfig{file=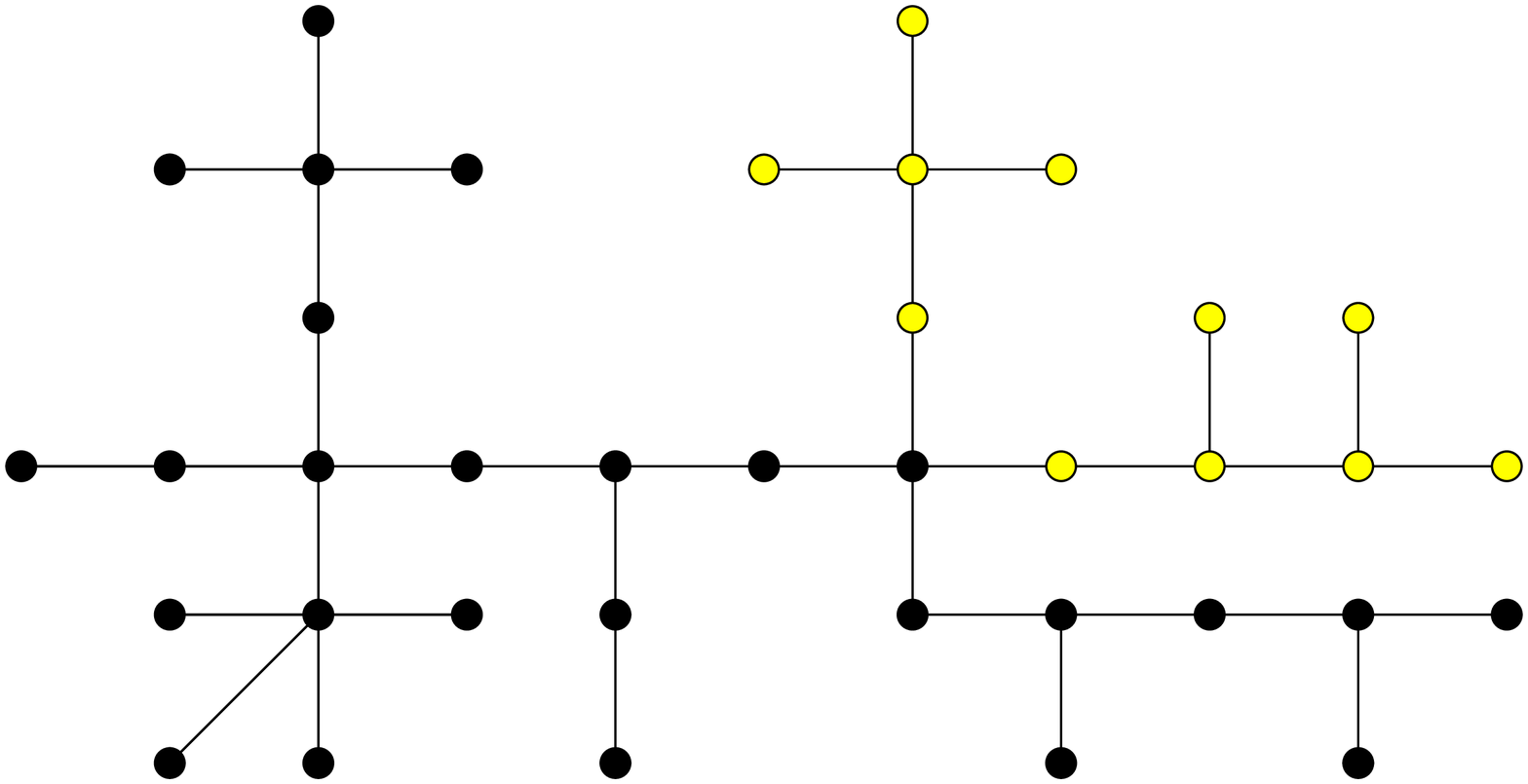,scale=0.25}\hspace*{30mm}
\epsfig{file= 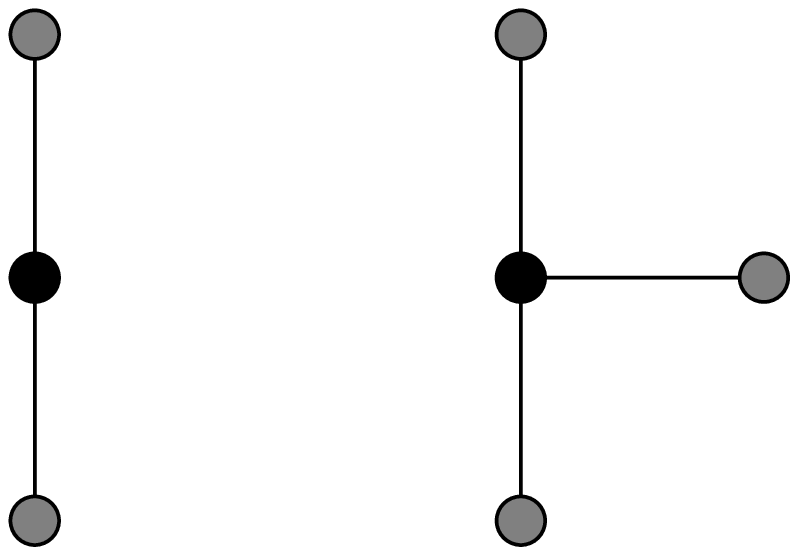,scale=0.4}
\end{center}
\caption{Eigenvector zero-nonzero patterns of a tree and corresponding skeleton forest}\label{figfullexampleeval2}
\end{figure}

In Figure \ref{figfullexampleeval2} an example of a tree $T$ with threefold eigenvalue $2$ is shown along
with its skeleton forest $\skel{T}{2}$. The black vertices of $T$ denote the vertices on which the respective
eigenvector vanishes. It can be clearly seen how the respective components of $T\setminus \nzerox{T}{x}$ correspond 
to a part of the skeleton. The black vertices in the skeleton correspond to the set $\nzeroc{T}{2}$.

\begin{env}{Lemma}\label{cpartnoev}
Let $T$ be a tree with eigenvalue $\lambda$.
Then $\comp{T\setminus \nzeroc{T}{\lambda}}$ can be partitioned into $\comp{T\setminus \nzero{T}{\lambda}}$
and a set of trees without eigenvalue $\lambda$.
\end{env}

\begin{proof}
Considered as a subgraph of $T$ every component of $T\setminus \nzero{T}{\lambda}$ is adjacent only to
vertices from $\nzeroc{T}{\lambda}$, but by definition does not contain such vertices.
So $\comp{T\setminus \nzero{T}{\lambda}} \subseteq \comp{T\setminus \nzeroc{T}{\lambda}}$.
By construction all elements of $\comp{T\setminus \nzero{T}{\lambda}}$ have eigenvalue $\lambda$.

Now let $C\in\comp{T\setminus \nzeroc{T}{\lambda}}\setminus\comp{T\setminus \nzero{T}{\lambda}}$.
All vertices of $C$ necessarily belong to the set $\nzero{T}{\lambda}$ so
that every eigenvector of $T$ for eigenvalue $\lambda$ must vanish on $C$.

Assume that there exists an eigenvector $y$ of $C$ for eigenvalue $\lambda$. 
Construct an eigenvector $z$ for eigenvalue $\lambda$ of $T$ as follows. Firstly, 
set $z\vert_C = y$.
Consider a vertex $w$ that is adjacent to $C$ in $T$ and let $v$ be the neighbor of $w$ in $C$.
Let $\nu$ be the value of $y$ on $v$.

Case $\nu=0$: 
Simply set $z$ to zero on the vertices of the particular component of $T\setminus C$ that
contains $w$.

Case $\nu\not=0$:
Clearly, $w\in\nzeroc{T}{\lambda}$ so that 
by construction, $w$ has a neighbor $u\not=v$ that belongs to a component of $T\setminus \nzero{T}{\lambda}$
(since $v$ does not). There exists an eigenvector $x$ of $T$ for eigenvalue $\lambda$
that vanishes on $C$ and $w$ but does not vanish on $u$. We may assume \obda
that $x$ has value $-\nu$ on $u$.
Let $T_u$ be the branch of $T$ connected to $w$ via
$u$. Let $T_F$ be the union of the branches connected to $w$ via the
neighbors of $w$ different from $u,v$. Note that $T_F$ is nonempty since $x$ must fulfil the
summation rule at vertex $w$. Set $z\vert_{T_F}=0$ and $z\vert_{T_u}=x\vert_{T_u}$. Now the
summation rule holds for $w$ and all the vertices of $T_u$ and $T_F$.

Apply the described procedure for every eligible vertex $w$. After that the values of $z$ are completely
determined. This yields a valid eigenvector for eigenvalue $\lambda$ of $T$ that
does not vanish on $C$, a contradiction.
\end{proof}

Combining Lemma \ref{cpartnoev} with Corollary \ref{compunion} and Lemma \ref{restrictev} we can derive
the following useful statement:

\begin{env}{Lemma}\label{restrictevlambdacomp}
Let $T$ be a tree and $x$ an eigenvector  for eigenvalue $\lambda$ of $T$.
Then for every $C\in\comp{T\setminus \nzeroc{T}{\lambda}}$ the restriction $x\vert_C$
either has only zero entries or only non-zero entries. In the latter case
it constitutes an eigenspace basis of the subgraph $C$ of $T$.
\end{env}

The next lemma is concerned with eigenvectors of skeletons.

\begin{env}{Lemma}\label{skelvecextend}
Let $T$ be a tree and $x$ an eigenvector for eigenvalue $\lambda$ of $T$.
Then every vector from the null space of $\skelx{T}{x}$ can be trivially extended
to a vector from the null space of $\skel{T}{\lambda}$.
\end{env}

\begin{proof}
Let $T'$ be the subgraph of $T$ induced by $\nzerocx{T}{x}$ and the vertices of
$T\setminus \nzerox{T}{x}$. Then $\nzero{T'}{\lambda}=\nzerocx{T}{x}$ and
$\skel{T'}{\lambda}=\skelx{T}{x}$.
Within $T$, vertices of $T'$ that are adjacent to vertices not in $T'$ 
necessarily belong to $\nzerocx{T}{x}$ for otherwise the summation rule would be violated.

Now consider a vector from the null space of $\skel{T'}{\lambda}$.
Extend it to  $\skel{T}{\lambda}$ by adding zero entries on the additional vertices.
Clearly, non-zero entries only occur on the vertices not in $\nzero{T'}{\lambda}$. 
But these have no neighbors belonging to $\skel{T}{\lambda}$ but not to $\skel{T'}{\lambda}$. So the
summation rule is trivially fulfilled for all vertices of $\skel{T}{\lambda}$.
\end{proof}


\begin{env}{Lemma}\label{metaskel}
Let $T$ be a tree with eigenvalue $\lambda$. Further let $S'=T/\comp{T\setminus \nzeroc{T}{\lambda}}$ and $S=\skel{T}{\lambda}$.

Then the skeleton $S$ forms an induced sub-forest of the tree $S'$ such that 
$S'\setminus V(S)$ contains no edges.
\end{env}

\begin{proof}
Assume that $S$ does not form an induced subgraph of $S'$. Then two vertices of $S$ are adjacent
in $S'$ but not already in $S$. Since these vertices by construction must lie in the same
component of $S$, the additional edge would create a cycle in $S'$, which is impossible. 
By construction the vertices of $S'\setminus V(S)$ are mutually non-adjacent in $S'$.
\end{proof}

Lemma \ref{metaskel} allows us to derive the notion of a {\em meta skeleton} in which the components
of the skeleton forest of a tree $T$ with eigenvalue $\lambda$ are joined by exactly those vertices contracted from the component trees
of $\comp{T\setminus \nzeroc{T}{\lambda}}$ that do not have eigenvalue $\lambda$ (cf.~Lemma \ref{cpartnoev}).


Next we explore the relation between eigenspace bases of trees and null space bases of the respective skeleton forests.

\begin{env}{Construction}\label{constrstraight}
Let $B=\{b_1,\ldots,b_r\}$ be an eigenspace basis for eigenvalue $\lambda$ of a given tree $T$.
Construct a basis $B'=\{b'_1,\ldots,b'_r\}$ of the same eigenspace as follows.
Let initially $b'_i=b_i$ for $i=1,\ldots,r$ and let $M=\emptyset$. There exists a component $C_1\not\in M$ of $T\setminus \nzero{T}{\lambda}$
such that $b_1\vert_{C_1}$ does not vanish. Owing to Lemma \ref{restrictevlambdacomp} we can subtract suitable
multiples of $b'_1$ from $b'_2,\ldots,b'_r$ such that $b'_i\vert_{C_1}=0$ for $i=2,\ldots,r$. Add $C_1$ to the set $M$.
Proceed iteratively for $b_j$, $j=2,\ldots,r$, by in turn finding a suitable $C_j\not\in M$ and 
establishing $b'_i\vert_{C_j}=0$ for $i=j+1,\ldots,r$.
\end{env}

The previous construction immediately gives rise to the following observation.

\begin{env}{Observation}\label{multlessthannumcomp}
Let $T$ be a tree and let $\lambda$ be an eigenvalue of $T$ with multiplicity $r\geq 1$. Then,
\[
\vert \comp{T\setminus \nzero{T}{\lambda}} \vert \geq r.
\]
\end{env}

We say that a set $\{x_1,\ldots,x_r\}$ of eigenvectors for eigenvalue $\lambda$ of a tree $T$ is {\em straight} if 
the components of $T\setminus \nzero{T}{\lambda}$ can be numbered $C_1,\ldots,C_s$ such that
for $j=1,\ldots,r$ we have $x_j\vert_{C_j}\not=0$ but $x_j\vert_{C_i}=0$ for $i=j+1,\ldots,r$.
Observation \ref{multlessthannumcomp} guarantees that $s\geq r$. 
Note that by Lemma \ref{restrictevlambdacomp} each condition $x_j\vert_{C_j}\not=0$ actually
means that $x_j$ vanishes on none of the vertices of $C_j$.
By Construction \ref{constrstraight} every tree eigenspace has a straight basis.

\begin{env}{Lemma}\label{straightindep}
Every straight set of tree eigenvectors is linearly independent.
\end{env}

\begin{proof}
Let $X=\{x_1,\ldots,x_r\}$ be a straight set of eigenvectors for eigenvalue $\lambda$ of a tree $T$.
Let $u_1,\ldots,u_r$ be vertices of $T$ such that they belong to distinct subgraphs of $T$ representing components of $T\setminus \nzero{T}{\lambda}$.
There exist enough such vertices because of Observation \ref{multlessthannumcomp}. Let $M$ be the matrix formed by taking 
$x_1,\ldots,x_r$ as columns and then retaining only those rows that correspond to the entries 
of the $x_i$ on the vertices $u_1,\ldots,u_r$. If $X$ is straight then with respect to a suitable numbering
of the $x_i$ the vertices $u_i$ can be selected such that $M$ is a lower diagonal matrix with non-zero diagonal entries.
Hence, the set $X$ is linearly independent.
\end{proof}

\begin{env}{Theorem}\label{maptreetoskelvec}
Let $T$ be a tree with eigenvalue $\lambda$ and corresponding eigenspace basis $B$. 
Then for every vector $b\in B$ there exists a vector $b'$ from
the null space of the skeleton $\skel{T}{\lambda}$ such that $b'$
is non-zero exactly on the vertices corresponding to the 
contracted subgraphs of $T$ on which its associated vector $b\in B$ does not vanish.
If $B$ is straight, then the vectors created from $B$ are linearly independent.
\end{env}

\begin{proof}
Let $b\in B$ and initialize $b'=0$. 
In the following let $C(v)$ denote the 
contracted subgraph corresponding to a vertex $v$ of $\skel{T}{\lambda}\setminus\nzeroc{T}{\lambda}$.
Moreover, if two vertices from $\skel{T}{\lambda}\setminus\nzeroc{T}{\lambda}$ have a common
neighbor in $\skel{T}{\lambda}$ (necessarily from $\nzeroc{T}{\lambda}$) they are called {\em brothers}.

For every component $S$ of $\skel{T}{\lambda}$ proceed as follows.
Fix a vertex $s$ of $S\setminus\nzeroc{T}{\lambda}$. If $b$ is non-zero on $C(s)$,
then $b'$ to $1$ on $s$. Consider $s$ as visited and all other vertices of $S\setminus\nzeroc{T}{\lambda}$ 
as unvisited. We now employ a tree search that starts at $s$ and iteratively corrects the values 
of $b'$ on the vertices of $S\setminus\nzeroc{T}{\lambda}$ such that, finally, $b'\vert_S$ belongs
to the null space of $S$ and assumes the desired zero-nonzero pattern. The search only visits
unvisited brothers of already visited vertices.

Let $v$ be a visited vertex of $S\setminus\nzeroc{T}{\lambda}$ that has unvisited brothers.
Mark all brothers of $v$ as visited once the steps described below have been carried out.
Let $W\subseteq\nzeroc{T}{\lambda}$ contain all vertices that are adjacent (in $\skel{T}{\lambda}$) both to $v$
and some unvisited brother of $v$. Now iterate over the vertices $w\in W$. 
Let $v_1,\ldots,v_r$ be all those unvisited brothers of $v$ that are adjacent to $w$ and 
for which $b$ does not vanish on $C(v_i)$. By construction, each vertex $v_i$ has exactly
one visited brother, namely $v$. 
Observe at this point that we have necessarily $r\geq 1$ if $b$ does not vanish on $C(v)$ because
else the summation rule would fail for $b$ on the vertex in $T$ that corresponds to $w$.
Hence, it is always possible to assign suitable non-zero values to the vertices $v_1,\ldots,v_r$ such that
$b'$ fulfils the summation rule for vertex $w$.

By construction and the definition of a skeleton it follows immediately that $b'$ is a valid eigenvector
from the null space of $\skel{T}{\lambda}$. Its zero-nonzero pattern is as claimed.
If $B$ is straight, then the set of vectors created from all the vectors of $B$ using the above procedure is straight as well.
Therefore it is linearly independent by Lemma \ref{straightindep}.
\end{proof}

Since every tree eigenspace has a straight basis we can immediately relate the dimensions of a tree eigenspace 
and the null space of the associated skeleton:

\begin{env}{Corollary}\label{multltnulldim}
Let $T$ be a tree with eigenvalue $\lambda$ of multiplicity $r\geq 1$. Let $s$ be the 
nullity of $\skel{T}{\lambda}$. Then $r\leq s$.
\end{env}


\begin{env}{Theorem}\label{mapskeltotreevec}
Let $T$ be a tree with eigenvalue $\lambda$ and let $B'$ be a basis
of the null space of its skeleton $\skel{T}{\lambda}$. 
Then for every vector $b'\in B'$ there exists an eigenvector $b$ of $T$ for
eigenvalue $\lambda$ such that $b$
is non-zero exactly on those subgraphs of $T$ that correspond to
vertices of $\skel{T}{\lambda}$ on which $b'$ does not vanish.
If $B'$ is straight, then the vectors created from $B'$ are linearly independent.
\end{env}

\begin{proof}
In view of Lemma \ref{restrictevlambdacomp} it is possible to use a technique similar to the one
used in the proof of Theorem \ref{maptreetoskelvec}, just in the opposite direction. 
\end{proof}

\begin{env}{Corollary}\label{multgtnulldim}
Let $T$ be a tree with eigenvalue $\lambda$ of multiplicity $r\geq 1$. Let $s$ be the 
nullity of $\skel{T}{\lambda}$. Then $r\geq s$.
\end{env}


By Corollaries \ref{multltnulldim} and \ref{multgtnulldim} we see that the multiplicity of the
eigenvalue $\lambda$ of a tree $T$ equals the nullity of the skeleton $\skel{T}{\lambda}$. 
It is well known that the nullity of a forest is closely linked to its matching properties.
We will exploit these ties with respect to skeletons.
But first let us generally relate maximum matchings of trees to eigenvectors of their null spaces.
Maximum matchings of trees can be quite elegantly obtained using specialized algorithms \cite{fricke96},\cite{fox1}.

\begin{env}{Lemma}\label{maymissnonadj}
Let $T$ be a tree.  Let $Q$ contain all vertices of $T$ such that each of them may be missed by some maximum matching of $T$.
Then the vertices of $Q$ are mutually non-adjacent in $T$.
\end{env}

\begin{proof} Given two given adjacent vertices $a,a'$, let $M,M'$ be two maximum matchings of $T$ such that $M$ does not cover $a$ and $M'$ 
does not cover $a'$.
Consider the two components that arise from the removal of the edge $aa'$ from $T$. Let $T_a$ denote the one that contains $a$ and $T_{a'}$
the one that contains $a'$. The matchings $M,M'$ each induce maximum matchings on both $T_{a}$ and $T_{a'}$. Construct a matching of $T$
by joining the matching induced by $M$ on $T_{a}$ and the one induced by $M'$ on $T_{a'}$. The resulting matching has
the same number of edges as $M,M'$ but can be extended by adding the edge $aa'$, a contradiction.
\end{proof}

\begin{env}{Theorem}\label{kernelbasis}
Let $T$ be a tree with edge set $E$. Let $K$ contain all vertices of $T$ that may be missed by maximum matchings of $T$.
Further, let $N$ contain all vertices of $T$ that are never missed by any maximum matching of $T$.
Consider a fixed maximum matching $M$ of $T$ and let $K_M\subseteq K$ be the vertices actually missed by the matching.

Then a simply structured null space basis of $T$ can be constructed as follows. Pick a vertex $v\in K_M$ and find the
subtree $S_v$ of $T$ formed by the union of all maximal paths that start at $v$ and alternatingly contain
edges from $E\setminus M$ and $M$, such that each edge in the path is incident to one vertex from $N$ and one from $K\setminus (K_M\setminus \{v\})$.
Assign weight $1$ to all vertices of $S_v$ whose distance to $v$ is divisible by four, assign weight $-1$ if the
distance is two modulo four, and assign zero to all other vertices of $T$.
\end{env}

\begin{proof}
Let us first verify that the summation rule holds on the tree $S_v$.
By construction, the vertices of $S_v$ receive a non-zero weight if and only if they belong to $K\setminus (K_M\setminus \{v\})$,
whereas the zero weight vertices of $S_v$ all belong to the set $N$. 
Let $w$ be a vertex of $S_v$. If $w\in K\setminus (K_M\setminus \{v\})$, then $w$ has only 
neighbors belonging to the set $N$, so that the summation rule is trivial to check.
If $w\in N$, then all its neighbors are from $K\setminus (K_M\setminus \{v\})$. However,
$w$ has necessarily degree $2$ in $S_v$ since otherwise either the given matching would not be valid
or a path using two consecutive edges from $E\setminus M$ would have been used in the creation of $S_v$.
These two neighbors of $w$ in $S_v$ have values $1$ and $-1$ so that the summation rule
holds for $w$.

In order to verify the summation rule on $T$ we only need to assert that no vertex $w$ of $S_v$ 
that belongs to $K\setminus (K_M\setminus \{v\})$ has a neighbor $x$ in $T$ that does not belong to $S_v$.
Assume to the contrary, that such vertices $w,x$ exist. Now $x$ would either belong to $K_M$, in which case an augmenting
path from $v$ to $x$ would exist in $T$ and therefore contradict the maximality of $M$. Or there
would exist an edge $xy\in M$ with $y\not=w$, contradicting the construction of $S_v$.

Linear independence of the constructed vectors is obvious by construction. That indeed a basis is formed
follows from the fact that the rank of a tree equals twice the number of edges in a 
maximum matching of the tree (see e.g.~\cite{bevisdom95}, \cite{fricke96}).
\end{proof}

\begin{env}{Corollary}\label{kerndimmissvanish}
Let $T$ be a tree. Then the set of vertices never missed by a maximum matching of $T$ is exactly
the set of vertices on which every vector from the null space of $T$ vanishes.
Moreover, the nullity of $T$ equals the number of vertices missed by a maximum matching of $T$. 
\end{env}

The second part of Corollary \ref{kerndimmissvanish} can also be found in \cite{spectraofgraphs} and \cite{fricke96}.

\begin{env}{Corollary}\label{nylenlike}
Let $T$ be a tree and let $R$ be the set of those vertices of $T$ on which the null space of $T$
does not completely vanish. Then the nullity of $T$ equals
the number of connected components of the subgraph of $T$ induced by the set $R$
minus the number of vertices of $T$ that are adjacent to $R$ but
not contained in it.
\end{env}

We will revisit Corollary \ref{nylenlike} later on in section \ref{sectreepattern}.

\begin{env}{Theorem}\label{skelmatch}
Let $T$ be a tree with eigenvalue $\lambda$. 

Then the set of vertices of the skeleton $\skel{T}{\lambda}$ that may be missed by a maximum matching of the skeleton
consists exactly of the vertices corresponding to the contracted components of $T\setminus \nzero{T}{\lambda}$.

The number of vertices of $\skel{T}{\lambda}$ that are missed by a maximum matching of the skeleton
equals the multiplicity of eigenvalue $\lambda$ of $T$.

The non-zero entries of a vector from the null space of $\skel{T}{\lambda}$ only occur on 
vertices that correspond to the contracted elements of $\comp{T\setminus \nzero{T}{\lambda}}$.
\end{env}

\begin{proof}
This follows from Corollary \ref{multltnulldim}, Corollary \ref{multgtnulldim},
Theorem \ref{kernelbasis} and Corollary \ref{kerndimmissvanish}.
\end{proof}

\begin{env}{Remark}\label{reconstrskel}
Given a tree, a $\lambda$-skeleton and a maximum matching of the skeleton, 
we can determine the multiplicity of $\lambda$ as an eigenvalue of $T$ by Theorem \ref{skelmatch}. 
Moreover, since Theorem \ref{kernelbasis} allows to construct a straight basis of the skeleton null space 
from a given maximum matching, a basis of the corresponding tree eigenspace can 
be obtained constructively by means of Theorem \ref{mapskeltotreevec}.

The only catch is that usually the skeleton vertices that may or may never be missed by a 
maximum matching are not known beforehand. Luckily, the so-called FOX algorithm presented 
in \cite{fox1}, \cite{fox2} allows to generate maximum matchings along with the
sets $K$,$N$ as required by Theorem \ref{kernelbasis}. Every maximum matching can be
obtained by a suitable run of the FOX algorithm.
\end{env}

Now we derive an interesting feature of a skeleton that becomes important once we 
want to characterize what forests may actually occur as skeletons (as we will see later in the
proof of Theorem \ref{blowupmetaskel}):

\begin{env}{Lemma}\label{skelnoeverymatchedge}
A skeleton $\skel{T}{\lambda}$ of a tree $T$ with eigenvalue $\lambda$ does not contain any edges
that belong to every maximum matching of the skeleton.
\end{env}

\begin{proof}
Theorem \ref{maptreetoskelvec},
Theorem \ref{kernelbasis} and the definition of $\skel{T}{\lambda}$
imply that a pair of skeleton vertices can only be adjacent if 
one of them is never missed by a maximum matching whereas the other one
may be missed by a maximum matching of the skeleton.
\end{proof}


Concluding this section, it is interesting to note that the skeleton construction cannot be arbitrarily iterated
in the sense that a skeleton is its own eigenvalue $0$ skeleton:

\begin{env}{Lemma}
Let $T$ be a tree with eigenvalue $\lambda$ and let $S=\skel{T}{\lambda}$ its $\lambda$-skeleton.
Then the skeleton $\skel{S}{0}$ equals $S$.
\end{env}

\begin{proof}
Clearly zero is an eigenvalue of the skeleton $S$. 
Every vertex of $S$ corresponds to a subgraph of $T\setminus \nzero{T}{\lambda}$ for which there exists an
eigenvector of $T$ for eigenvalue $\lambda$ that is non-zero on all vertices of the subgraph.
So by Theorem \ref{maptreetoskelvec} for every vertex of $S$ that corresponds to
an element of $\comp{T\setminus \nzero{T}{\lambda}}$ there exists a vector from the null space of $S$ 
that does not vanish on this vertex.

Consequently, the partition of the components of $S\setminus\nzeroc{S}{0}$ according to Lemma \ref{cpartnoev}
does not contain any trees without eigenvalue zero so that $\nzeroc{S}{0}= \nzero{S}{0}$.
On the other hand, a component of $S\setminus\nzeroc{S}{0}$ that has a null space eigenvector without zero entries 
must necessarily contain only a single vertex (cf.~Theorem \ref{kernelbasis}).
So every vertex of $S$ is associated with exactly one vertex in the skeleton $\skel{S}{0}$
and no factual contraction of subgraphs of $S$ happens when forming $\skel{S}{0}$.
\end{proof}

\subsection{Tree eigenvector composition} 

In this section a composition technique is outlined that allows to "blow up" a given potential skeleton
to a tree whose skeleton is indeed the initial graph. The bottom line is that we can not only use
the maximum matching properties of skeletons to describe eigenvectors of a tree but also conversely
construct a tree with predetermined eigenvector properties by generating it from a suitable skeleton.

Following the lines of Lemma \ref{metaskel}, Theorem \ref{skelmatch} and Lemma \ref{skelnoeverymatchedge} we introduce the following definition.
We call a tree a {\em meta skeleton} if its vertex set contains an independent set $X$
such that each subtree $C\in\comp{T\setminus X}$ does not have a perfect matching and only 
vertices of $C$ that are contained in every maximum matching of $C$ are adjacent to vertices of $X$ in $T$.
Moreover, it is required that no edge of $T\setminus X$ is contained in every maximum matching of $T\setminus X$.
The set $X$ is called an admissible {\em non-eigenvalue set} of the meta skeleton (cf.~Lemma \ref{cpartnoev}).

\begin{env}{Construction}\label{constrblowupmetaskel}
Let $S'$ be a meta skeleton with non-eigenvalue set $X$. Choose a number $\lambda\in\RR$ and construct a tree $T$ as follows.
Substitute each vertex of $X$ with a tree without eigenvalue $\lambda$. For every component of $S'\setminus X$ replace each
of its vertices that may be missed by a maximum matching of the component with
a tree that has eigenvalue $\lambda$ and a corresponding eigenvector without zero entries.
Whenever a vertex of an adjacent pair of vertices is substituted with a graph, a single arbitrary vertex of the
substituted graph is chosen to become connected to the other vertex of the pair.
\end{env}

\begin{env}{Theorem}\label{blowupmetaskel}
If Construction \ref{constrblowupmetaskel} succeeds for a given triplet $(S',X,\lambda)$, then the generated tree $T$
has eigenvalue $\lambda$. Its multiplicity equals the number of vertices missed in a maximum matching of $S'\setminus X$.
Moreover, $\skel{T}{\lambda}=S'\setminus X$ and $S'=T/\comp{T\setminus \nzeroc{T}{\lambda}}$.
\end{env}

\begin{proof}
We show that the skeleton and meta skeleton are as claimed. Then the rest of the theorem follows 
from Theorems \ref{maptreetoskelvec}, \ref{mapskeltotreevec}, \ref{kernelbasis}, \ref{skelmatch} as 
outlined in Remark \ref{reconstrskel}.

By the definition of a meta skeleton each component of $S'\setminus X$ has eigenvalue zero so that
we can use Theorem \ref{kernelbasis} to obtain a null space basis for each such component.
Let $K$ consist of all vertices of $S'$ that may be missed by maximum matchings of
their respective components of $S'\setminus X$.
With the same technique that was used in the proof of Theorem \ref{nzerocompcontain} we
can actually find a vector $y$ from the null space of the forest $\bigcup\comp{S'\setminus X}$
that is nonzero on every vertex of $K$. Let $T'$ be the subgraph of $T$ that is the
blown up counterpart of the subgraph $\bigcup\comp{S'\setminus X}$ of $S'$.
Using a technique similar to the one used in the 
proof of Theorem \ref{maptreetoskelvec} we can now employ the zero-nonzero pattern of $y$ 
to construct an eigenvector $x'$ for eigenvalue $\lambda$ of $T'$ since exactly the
vertices of $S'$ on which $y$ is non-zero have been blown up to suitable trees.
Use zero entries to trivially extend $x'$ to a vector $x$ on $T$. 
Then Corollary \ref{kerndimmissvanish} and the meta skeleton definition imply that $x$ is an eigenvector
of $T$ for eigenvalue $\lambda$ since the summation rule clearly also holds for the vertices 
of the subtrees blown up from the elements of $X$. 

We will see in a moment that the vector $x$ has been chosen such that it is zero exactly on the
vertices $\nzero{T}{\lambda}$, i.e.~$\nzerox{T}{x}=\nzero{T}{\lambda}$. This 
fact substantially eases the determination of the skeleton.

Let $N$ be the set of vertices of $S'\setminus X$ that are covered by every maximum matching of $S'\setminus X$. 
It follows from the definition of the meta skeleton and Lemma \ref{skelnoeverymatchedge}
that every vertex of $N$ is adjacent to some vertex of $S'\setminus X$ that may be missed by a maximum
matching of that graph. Since by Construction \ref{constrblowupmetaskel}
we can consider $N$ also as a subset of the vertices of $T$ it follows immediately 
that $N\subseteq\nzeroc{T}{\lambda}$.

Now assume that there exists an eigenvector $z$ for eigenvalue $\lambda$ of $T$ that is non-zero on a
subtree $T_v$ of $T$ blown up from a vertex of $v\in X$. Since by construction all outer neighbors of $T_v$
in $T$ belong to $N$ the restriction $z\vert_{T_v}$ must be a valid eigenvector for eigenvalue $\lambda$ of $T_v$.
But this is impossible by the choice of $T_v$. 
Hence, we have $N=\nzeroc{T}{\lambda}$, $\nzerox{T}{x}=\nzero{T}{\lambda}$ and $\comp{T\setminus \nzerox{T}{x}} = \comp{T\setminus \nzero{T}{\lambda}}$.

Moreover, if $T_X$ denotes the sub-forest of $T$ that is the union of all graphs blown up from the vertices of $X$,
then $\nzero{T}{\lambda}=V(T_X)\cup N$. 
So $S'\setminus X$ is indeed the skeleton of $T$ and $S'=T/\comp{T\setminus \nzeroc{T}{\lambda}}$.
\end{proof}

\begin{figure}[ht]
\begin{center}
\epsfig{file=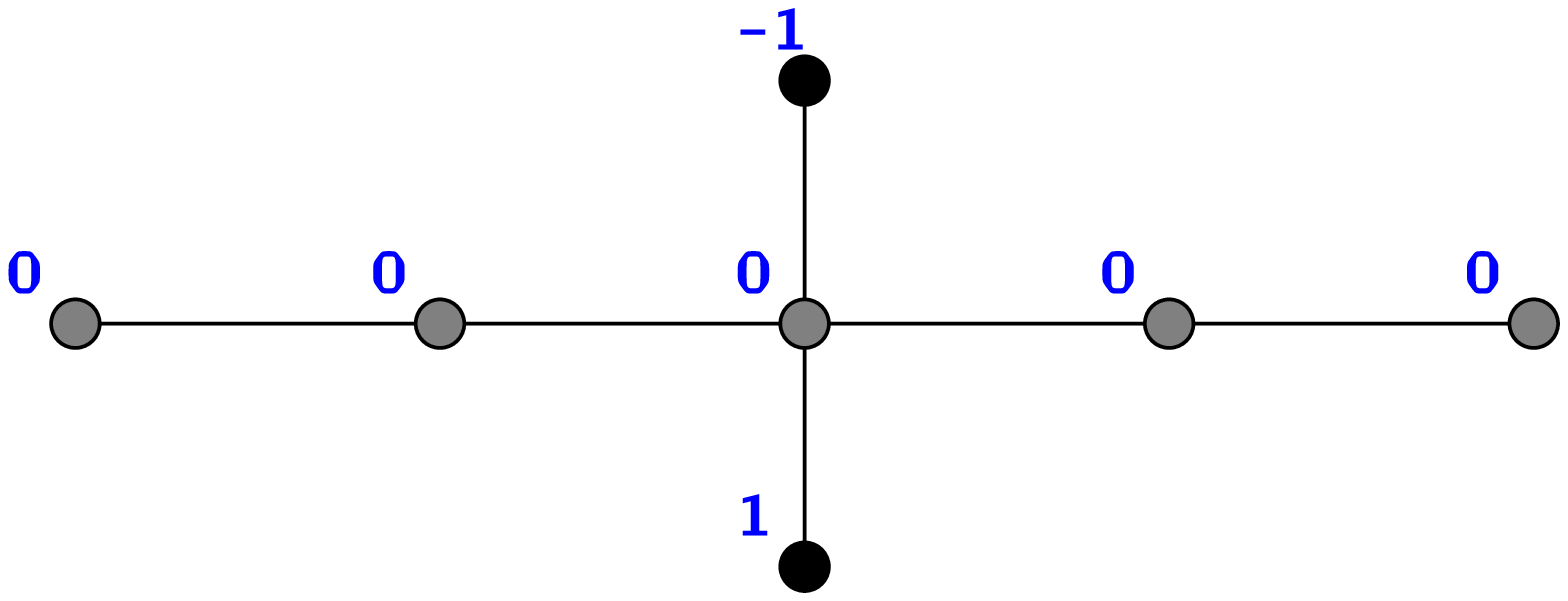,scale=0.35}
\end{center}
\caption{Malformed skeleton example} \label{figmalformedskel}
\end{figure}

\begin{env}{Remark}
The skeleton property stated in Lemma \ref{skelnoeverymatchedge} is decisive for choosing proper
forests to be blown up. Otherwise, even though valid eigenvectors can be constructed for the blown up graph, 
it cannot be guaranteed that a proper eigenspace basis is obtained because the skeleton of the blown up graph
may in fact not be the graph we expanded. See Figure \ref{figmalformedskel} for a malformed skeleton with
one-dimensional null space that can be blown up to the graph with eigenvalue $1$ shown in Figure \ref{figsameskel}.
However, its eigenvalue $1$ has multiplicity $3$.
\end{env}

\section{Applications}

\subsection{Simply Structured Tree Eigenspace Bases} 

In this section we give a characterization of trees for which eigenspace bases with entries only from $\{0,1,-1\}$ exist.
Let us call a vector space basis {\em simply structured} if it only contains vectors with entries from the set $\{0,1,-1\}$.

It is easy to see that the only feasible eigenvalues that allow the construction of such bases are $0,1,-1$.
Simply apply the summation rule to a leaf with non-zero value. It follows by a straightforward argument that
such a leaf exists for every non-null eigenvector.

It has already been independently shown in \cite{akbariali06},\cite{fox1} that every tree has such a simply structured basis
for eigenvalue $0$.  We complete the characterization by investigating the other two possible eigenvalues.
To this purpose we make use of the concept of decomposing trees by the zero entries of their eigenvectors that
was presented earlier.
Since trees are bipartite it now suffices to restrict further investigations to the eigenvalue $1$.
Given an eigenspace basis for eigenvalue $1$ an eigenspace basis for eigenvalue $-1$ is readily 
obtained by negating the signs of all vector entries corresponding to the vertices of one part of the
vertex bipartition.

\begin{figure}[ht]
\begin{center}
\epsfig{file=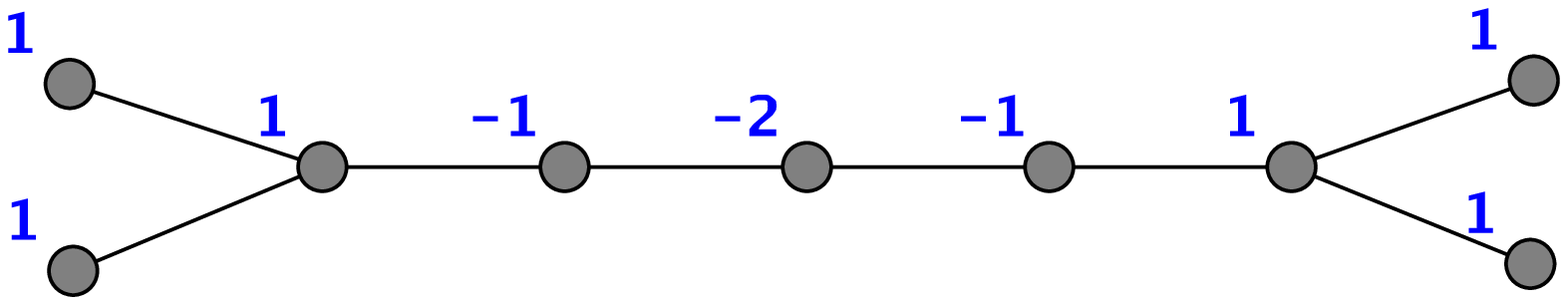,width=7cm}
\end{center}
\caption{Graph without $\{0,1,-1\}$ eigenvector for eigenvalue $1$} \label{fignonsimpleexample}
\end{figure}

Examples for eigenvectors for eigenvalue $1$ that cannot be scaled to $\{0,1,-1\}$ entries can
be found quite easily. See for example Figure \ref{fignonsimpleexample}, where the claim follows by Lemma \ref{lemfiedler}. 
In the following we will therefore
attempt to characterize those trees that have a simply structured eigenspace basis for eigenvalue $1$.

With respect to the existence of such trees we can make the following positive assertion.

\begin{figure}[ht]
\begin{center}
\epsfig{file=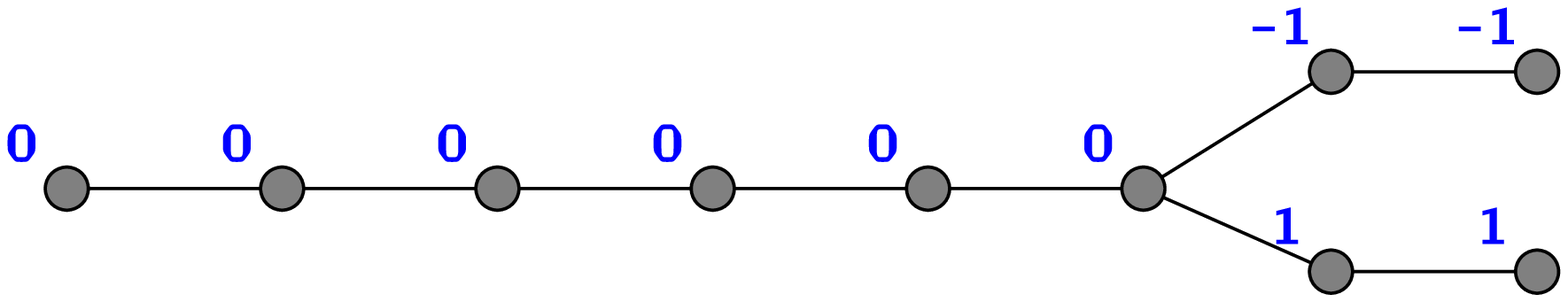,width=7cm}
\end{center}
\caption{Construction of trees with $\{0,1,-1\}$ eigenvector for eigenvalue $1$} \label{figsimpleconstrexample}
\end{figure}

\begin{env}{Observation}
For every integer $n\geq 5$ there exists a tree with $n$ vertices that has a
$\{0, 1,-1\}$ eigenvector for eigenvalue $1$. Simply take a path $P_5$ and attach
a pendant path $P_{n-5}$ to its center vertex, see Figure \ref{figsimpleconstrexample}.
In effect, further analysis of the construction shows that these trees 
admit simply structured eigenspace bases for eigenvalue $1$.
\end{env}

Assume that a tree with a simply structured eigenspace basis for eigenvalue $1$ is decomposed according to the
always-zero entries. Clearly, each such generated component has a single eigenvalue $1$ and a corresponding 
eigenvector without zero entries, namely a $\{1, -1\}$ vector. Since such eigenvectors are the building blocks
for the composition of trees with simply structured bases for eigenvalue $1$ we now direct our attention to them.
It turns out that trees with $\{1, -1\}$ eigenvector for single eigenvalue $1$ can be characterized in a very elegant way.

\subsubsection{$\{1,-1\}$ Eigenvectors for Eigenvalue $1$} 

\begin{env}{Observation}\label{leafandneigh}
Let $x$ be an eigenvector for eigenvalue $1$ of a given tree $T$. Then the value of $x$ on a
leaf equals the value on its unique neighbor.
\end{env}

\begin{figure}[ht]
\begin{center}
\epsfig{file=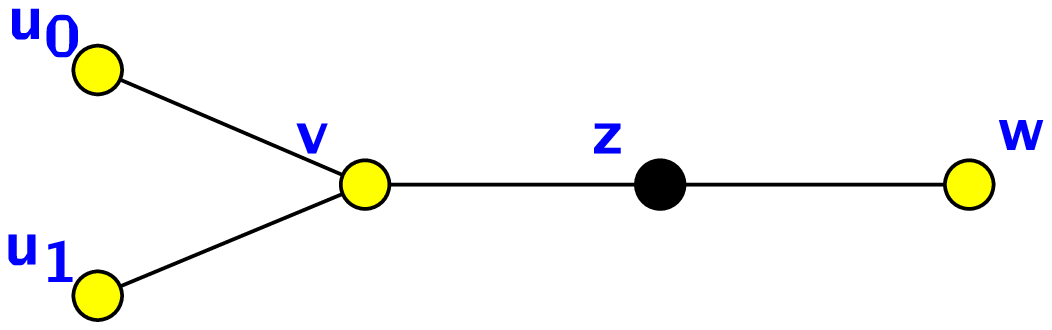,width=5cm}\hspace*{2cm}\epsfig{file=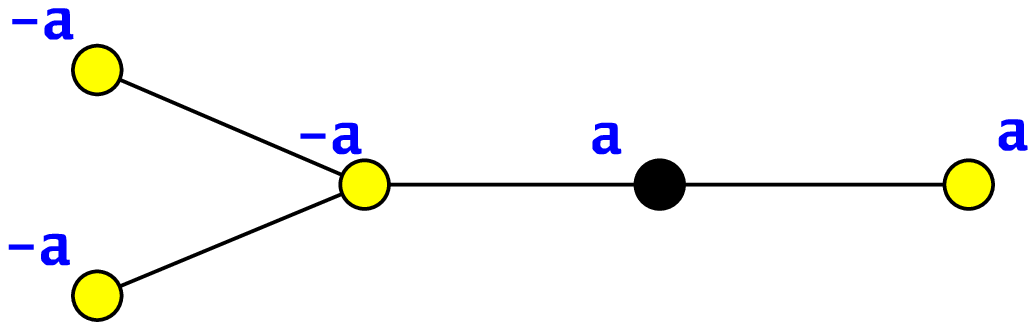,width=5cm}
\end{center}
\caption{Reduction subgraph and weights for $\{1,-1\}$ eigenvectors} \label{figeval1extension}
\end{figure}

\begin{env}{Theorem}\label{ew1reductionchar}
A tree has a $\{1,-1\}$ eigenvector for eigenvalue $1$ if and only if the tree can be reduced 
to a $K_2$ graph by repeatedly selecting a subgraph as in Figure \ref{figeval1extension} (where
the vertices $u_0,u_1,w$ must be leaves in the current reduced graph) and removing all its vertices
except $v$ from the current reduced graph.
\end{env}

\begin{proof}
Let $T$ be a tree with $\{1,-1\}$ eigenvector $x$ for eigenvalue $1$. Clearly, $T$ must have at least two
vertices. If $T$ is a complete graph $K_2$ there is nothing to show. So we may assume that $T$ has at least three vertices.

Let $u_0$ be a leaf of $T$ that is the end vertex of a diameter path and $v$ its only neighbor. Among those neighbors of $v$ different
from $u_0$ let $z$ be the uniquely determined vertex that is closest to the center of $T$. Let $u_1,\ldots,u_r$ be
the other neighbors of $v$ besides $u_0$ and $z$. Since $u_0$ has maximum eccentricity the vertices $u_1,\ldots,u_r$
must also be leaves of $T$.

We may assume that $v$ is not the sole center vertex of $T$. Otherwise $T$ would be a star graph $K_{1,r+2}$, which
does not have eigenvalue $1$. Let \obda~$x$ have value $1$ on $u_0$.
Then by Observation \ref{leafandneigh}, $x$ assumes the same value also on the vertices $u_1,\ldots,u_r,v$.
The summation rule for vertex $v$ requires a negative value of $x$ on $z$. Therefore, $r=1$ and $x$ has
value $-1$ on $z$. 

We now claim that $z$ is adjacent to a leaf with value $-1$. By the summation rule there exist at least
two neighbors of $z$ on which $x$ assumes the value $-1$. Among these neighbors there exists at least
one vertex $w$ such that the branch adjacent to $z$ via the edge $wz$ does not contain any center vertices
of $T$. Assume that $w$ is not a leaf of $T$. Then by the summation rule $w$ would have at least one neighbor
$w'$ with value $1$. Again by the summation rule $w'$ would have at least one neighbor $w''$ with value $1$.
But by our assumption about the location of the center of $T$ the eccentricity of $w''$ is clearly 
greater than that of $u_0$, a contradiction.

Remove the vertices $u_0,u_1,v,w$ from $T$. Clearly, $T$ remains a tree. Moreover, the summation rule remains valid for all remaining 
vertices, in particular for $z$. Since $z$ has at least one remaining neighbor it follows that $T$
has at least two vertices. We can therefore iterate the reduction step until a graph $K_2$ has been obtained.
The reduction procedure can also be applied for every subgraph of $T$ isomorphic to the one in Figure \ref{figeval1extension}
if only $u_0,u_1,w$ are leaves. The maximum eccentricity criterion only asserts the existence of a subtree suitable for reduction.

Conversely, assume that a tree can be decomposed in the described manner. Then we can assemble it from $K_2$ by
iteratively selecting a vertex $z$ and adding vertices $u_0,u_1,y,w$ according to Figure \ref{figeval1extension}.
The all ones vector forms an eigenspace basis for eigenvalue $1$ of the graph $K_2$. 
After the addition of the vertices  $u_0,u_1,y,w$ we can uniquely augment the previous eigenvector to become
a $\{1,-1\}$ eigenvector for eigenvalue $1$ of the extended graph. The values on the newly added vertices depend
only on the existing eigenvector value on $z$, cf.~Figure \ref{figeval1extension}. Iterating this
argument we find that $T$ has a $\{1,-1\}$ eigenvector for eigenvalue $1$.
\end{proof}

\begin{env}{Corollary}
There exists a tree with $n$ vertices that has a $\{1,-1\}$ eigenvector for eigenvalue $1$ if and only if  
$n\equiv 2 \mod 4$.
\end{env}

In the following, let $\cclass$ denote the class of all trees with $\{1,-1\}$ eigenvector for eigenvalue $1$.
Note that if a tree has a $\{1,-1\}$ eigenvector for eigenvalue $1$, then the eigenvalue $1$ has necessarily multiplicity one.

\subsubsection{$\{0,1,-1\}$ Eigenvectors for Eigenvalue $1$} 

Having investigated trees with $\{1,-1\}$ eigenvectors it is now straightforward
to achieve a characterization of trees with simply structured eigenspace bases for eigenvalue $1$:

\begin{env}{Theorem}\label{simplebasischar}
Let $T$ be a tree with eigenvalue $1$. Then there exists a simply structured basis for the
corresponding eigenspace if and only if $C\in\cclass$ for every component $C\in\comp{T\setminus \nzero{T}{1}}$.
\end{env}

\begin{proof}
Necessity follows from Lemma \ref{restrictevlambdacomp}. 
Sufficiency, on the other hand, follows from Remark \ref{reconstrskel} and Theorem \ref{blowupmetaskel}.
Just consider the reconstruction of (linearly independent) eigenvectors of $T$ from the zero-nonzero patterns of a straight
null space basis of its skeleton forest. Simply assign valid $\{1,-1\}$ eigenvectors to all contracted
subgraphs of $T$ where the chosen skeleton null space eigenvector is nonzero on the corresponding
skeleton vertices. A valid eigenvector is obtained by establishing the summation rule for all 
vertices of $\nzeroc{T}{1}$. This can be achieved by conducting a breadth first search from a 
fixed nonzero skeleton vertex $v$. Each time a vertex of $\nzeroc{T}{1}$ is visited
the summation rule for its partner vertex in $T$ is enforced by suitably multiplying the values on the branches 
leading away from $v$. Since the branches have only values from the set $\{0,1,-1\}$ the only factors
that are needed are $1$ and $-1$ so that finally a $\{0,1,-1\}$ eigenvector is created.
\end{proof}

In theory, Theorem \ref{simplebasischar} provides us with a completely structural characterization of all
trees whose eigenspace for eigenvalue $1$ admits a simply structured basis. The class $\cclass$ is
characterized by a reduction property and the set $\nzero{T}{1}$ is independent of the choice of a particular eigenspace
basis so that essentially it is an intrinsic structural property of a tree as well.

From a practical point of view, however, there is always an algebraic aspect.
In order to check if a tree $T$ has a simply structured eigenspace one would start by computing
an arbitrary eigenspace basis for eigenvalue $1$ and then try to reduce the components of 
$T\setminus \nzero{T}{1}$. Conversely, the trees with simply structured eigenspace bases for eigenvalue $1$ can be generated using
Construction \ref{constrblowupmetaskel} by using only graphs from $\cclass$ for the blown-up trees with eigenvalue $1$. 
But for blowing up the vertices of the non-eigenvalue set, trees without eigenvalue $1$ are used. So far a non-algebraic
characterization of such trees is unknown. It is even doubtful if such a characterization exists since it
is not difficult to show that every given tree can be extended to a tree with single eigenvalue $1$ and  a
corresponding eigenvector without zero entries. So it seems hard to tell the difference between trees
that have eigenvalue $1$ and trees that haven't. All in all, the desired characterization of trees
with simply structured eigenspace bases has been achieved.

We conclude this section with a construction that allows to derive a simply structured tree eigenspace basis
from a given initial eigenspace basis.

\begin{env}{Construction}
Let $T$ be a tree with eigenvalue $1$ and $B$ a corresponding eigenspace basis. Then a simply structured
eigenspace basis for this eigenvalue of $T$ can be obtained as follows:
\begin{enumerate}
\item Use $B$ to determine $\comp{T\setminus \nzero{T}{1}}$ and the skeleton $\skel{T}{\lambda}$.
\item Reduce every component of $T\setminus \nzero{T}{1}$ according to Theorem \ref{ew1reductionchar} and simultaneously
  determine $\{1,-1\}$ component eigenvectors.
\item Determine a maximum matching of $\skel{T}{\lambda}$, e.g.~using one of the algorithms presented in 
  \cite{bevisdom95}, \cite{fricke96} or \cite{fox1}.
\item Construct a straight null space basis $B'$ of $\skel{T}{\lambda}$.
\item Map the vectors of $B'$ to a set of vectors on $T$ by matching their zero-nonzero patterns;
  to the subgraphs corresponding to nonzero skeleton vertices the respective already computed 
  $\{1,-1\}$ component eigenvectors are assigned, all other vertices are assigned zero values.
\item For every constructed vector use a breadth first search approach to correct the signs of branches
such that for every zero value vertex adjacent to a non-zero vertex the summation rule holds.
\end{enumerate}
\end{env}

The initial basis $B$ can be obtained by traditional Gaussian elimination, but there exists an
algorithm that allows to compute the vectors of $B$
on the tree $T$ itself \cite{johnschild96}.

\subsection{Tree pattern matrices}\label{sectreepattern}

Let $M$ be a real $n\times n$ matrix. We associate with it a (directed) graph $\patterngr{M}$ with vertices
$v_1,\ldots,v_n$ such that there is an edge from $v_i$ to $v_j$ if and only if $M$ has a non-zero entry at
position $(i,j)$. If $\patterngr{M}$ is a tree, then we call $M$ a {\em tree pattern matrix}.
Let $\support{M}{\lambda}$ denote the set of vertices of $\patterngr{M}$ on which the eigenspace 
for eigenvalue $\lambda$ of $M$ does not entirely vanish. We call this set the  {\em support}
of $M$ with respect to $\lambda$. We define the support $\support{G}{\lambda}$ of a graph $G$
as the support of its adjacency matrix. Note that it is easy to find examples such that
 $\support{M}{\lambda}$ and $\support{\patterngr{M}}{\lambda}$
are different. 

We can extend Corollary \ref{nylenlike} to the following result which has 
already been published in \cite{nylen98} but proven differently:

\begin{env}{Theorem}\label{theoremnylen}
Let $M$ be an $n\times n$ tree pattern matrix. Then the dimension of the null space of $M$ equals
the number of connected components of the subgraph of $\patterngr{M}$ induced by $\support{M}{0}$
minus the number of vertices of $\patterngr{M}$ that are adjacent to $\support{M}{0}$ but
do not belong to this set.
\end{env}

\begin{proof}
Let $M$ be a tree pattern matrix and let $A$ be the adjacency matrix of $\patterngr{M}$. 
Theorem \ref{kernelbasis} states that $\support{A}{0}$ forms an independent vertex set in $\patterngr{M}$. 
Given a vector $v$ from the null space of $A$ we can transform it to a vector $v'$ from the
null space of $M$ having the same zero-nonzero pattern as follows. Assign $v$ to the vertices of  $\patterngr{M}$.
Conduct a breadth first search on $\patterngr{M}$ from a fixed vertex $s$ and enforce new summation rules. 
To be precise, for every vertex $z$ (as traversed by the breadth first search) it is possible
to multiply each of its adjacent branches leading away from $s$ with a nonzero factor such that 
the summation rule given by the line of $M$ that corresponds to $z$ holds. From a straight 
basis of the null space of $A$ we can thus obtain a straight linearly independent set of vertices from
the null space of $M$. A similar conversion can be employed for the opposite direction.
Therefore, $\support{M}{0}=\support{A}{0}=\support{\patterngr{M}}{0}$. 
Now the result follows by Corollary \ref{nylenlike}.
\end{proof}

In fact, the results from the previous sections allow us to generalize even further. 
We quoted Lemma \ref{lemfiedler} only as a special case of what is actually proven in \cite{fiedler75}. 
It has been shown that every eigenvector of a tree pattern matrix necessarily belongs to an eigenvalue with multiplicity one
if it has no zero entries. Moreover, for the application of the summation rule none of the proofs
given in section \ref{chmain} explicitly relied on the fact that it was induced by the adjacency matrix
of the tree. Every row of a tree pattern matrix $M$ induces a particular summation rule
for the associated vertex of $\patterngr{M}$. The only difference to the summation rule used for the adjacency matrix
is that for every vertex certain non-zero factors are applied to the weights of the neighbors before adding them up.
Consequently, we can generalize the entire theory presented in section \ref{chmain} to cover eigenvectors of
tree pattern matrices. In particular we obtain the following generalization of Theorem \ref{theoremnylen}:

\begin{env}{Theorem}\label{theoremnylenext}
Let $M$ be an $n\times n$ tree pattern matrix with eigenvalue $\lambda$. Then the dimension of the eigenspace of $M$ for eigenvalue
$\lambda$ equals the number of connected components of the subgraph of $\patterngr{M}$ induced by $\support{M}{\lambda}$
minus the number of vertices of $\patterngr{M}$ that are adjacent to $\support{M}{\lambda}$ but
do not belong to this set.
\end{env}

One other noteworthy generalization is that eigenspace dimensions of
tree patterned matrices are determined by sizes of maximum matchings of the respective associated 
skeletons.




\footnotesize
\bibliographystyle{acm}            
\bibliography{treedecompos_v2}  

\end{document}